\documentclass[numreferences]{kluwer}

\usepackage{amsmath}
\usepackage{xspace}

\newcommand{\rset}{\ensuremath{\mathbb{R}}\xspace}
\newcommand{\nset}{\ensuremath{\mathbb{N}}\xspace}
\newcommand{\rmax}{\ensuremath{\rset_{\max}}\xspace}
\newcommand{\rmin}{\ensuremath{\rset_{\min}}\xspace}

\newcommand{\Mat}{\mathrm{Mat}}

\newcommand{\0}{\ensuremath{\mathbf{0}}\xspace}
\newcommand{\1}{\ensuremath{\mathbf{1}}\xspace}

\newcommand{\cle}{\preccurlyeq}
\newcommand{\clt}{\prec}
\newcommand{\cgt}{\succ}

\newcommand{\oodot}{\ensuremath{\mathbin{\overline{\odot}}}\xspace}
\newcommand{\ooplus}{\ensuremath{\mathbin{\overline{\oplus}}}\xspace}
\newcommand{\oostar}{\ensuremath{\mathbin{\overline{\star}}}\xspace}
\newcommand{\bigooplus}{\mathop{\overline{\bigoplus}}}

\newcommand{\aint}{\mathbf a}

\newcommand{\A}{\mathbf A}
\newcommand{\lA}{\underline{\A}}
\newcommand{\uA}{\overline{\A}}

\newcommand{\B}{\mathbf B}
\newcommand{\lB}{\underline{\B}}
\newcommand{\uB}{\overline{\B}}

\newcommand{\lamint}{\boldsymbol{\lambda}}
\newcommand{\llambda}{\underline{\lamint}}
\newcommand{\ulambda}{\overline{\lamint}}

\newcommand{\tint}{\mathbf t}
\newcommand{\lt}{\underline{\tint}}
\newcommand{\ut}{\overline{\tint}}

\newcommand{\V}{\mathbf V}

\newcommand{\x}{\mathbf x}
\newcommand{\lx}{\underline{\x}}
\newcommand{\ux}{\overline{\x}}

\newcommand{\X}{\mathbf X}

\newcommand{\y}{\mathbf y}
\newcommand{\ly}{\underline{\y}}
\newcommand{\uy}{\overline{\y}}

\newcommand{\z}{\mathbf z}
\newcommand{\lz}{\underline{\z}}
\newcommand{\uz}{\overline{\z}}

\newcommand{\Z}{\mathbf{Z}}

\newtheorem{thm}{Theorem}
\newtheorem{prop}{Proposition}
\newtheorem{lm}{Lemma}

\newcounter{ExampleCount}[section]
\renewcommand{\theExampleCount}{\thesection.\arabic{ExampleCount}}
\newcommand{\example}[1]{%
\par%
\refstepcounter{ExampleCount}%
\textsc{Example \theExampleCount. #1{}}}

\newcounter{RemarkCount}[section]
\renewcommand{\theRemarkCount}{\thesection.\arabic{RemarkCount}}
\newcommand{\remark}{%
\par%
\refstepcounter{RemarkCount}%
\textsc{Remark \theRemarkCount. }}

\begin{document}

\begin{article}

\begin{opening}

\title{Idempotent Interval Analysis and Optimization Problems%
\thanks{The work was supported by the RFBR grant No.\ 99--01--01198 and the
Erwin Schr\"odinger Institute for Mathematical Physics.
\\ Submitted to \textit{Reliable Computing}.}}

\author{G.\ L.\ \surname{Litvinov}\email{litvinov@islc.msk.su}}
\institute{International Sophus Lie Centre}

\author{A.\ N.\ \surname{Sobolevski\u\i}\email{ansobol@idempan.phys.msu.su}}
\institute{M. V. Lomonosov Moscow State University}

\runningtitle{Idempotent Interval Analysis and Optimization Problems}
\runningauthor{G.\ L.\ Litvinov, A.\ N.\ Sobolevski\u{\i}}

\begin{ao}
\par\noindent
International Sophus Lie Centre\\
Nagornaya, 27--4--72\\
Moscow 113186 Russia\\
E-mail: litvinov@islc.msk.su
\end{ao}

\begin{abstract}
Many problems in optimization theory are strongly nonlinear in the
traditional sense but possess a hidden linear structure over suitable
idempotent semirings. After an overview of `Idempotent Mathematics' with an
emphasis on matrix theory, interval analysis over idempotent semirings is
developed. The theory is applied to construction of exact interval
solutions to the interval discrete stationary Bellman equation. Solution
of an interval system is typically $NP$-hard in the traditional interval
linear algebra; in the idempotent case it is polynomial. A generalization
to the case of positive semirings is outlined.
\end{abstract}


\classification{MSC codes}{65G10, 16Y60, 06F05, 08A70, 65K10}

\end{opening}

\section*{Introduction}%
\label{s:intro}

Many problems in the optimization theory and other fields of mathematics
are nonlinear in the traditional sense but appear to be linear over
semirings with idempotent addition.\footnote{One of the most important
examples of an idempotent semiring is the set $\rset_{\max} = \rset \cup
\{-\infty\}$ with operations $\oplus = \max$, $\odot = +$ (see
example~\ref{e:rmax} in subsection~\ref{s:semiringexamples}), sometimes
called the \emph{max-plus algebra}.} This approach is developed
systematically as \emph{Idempotent Analysis} or \emph{Idempotent
Mathematics} (see, e.g., \cite{Maslov}--\cite{BCOQ}). In this paper we
present an idempotent version of Interval Analysis (its classical version
is presented, e.g., in~\cite{Kearfott}--\cite{Moore}) and discuss
applications of the idempotent matrix algebra to discrete optimization
theory.

The idempotent interval analysis appears to be best suited for treating
problems with order-preserving transformations of input data. It gives
exact interval solutions to optimization problems with interval
uncertainties without any conditions of smallness on uncertainty intervals.
Solution of an interval system is typically $NP$-hard in the traditional
interval linear algebra; in the idempotent case it is polynomial. The
idempotent interval analysis is particularly effective in problems that are
strongly nonlinear in the traditional sense but possess a hidden linear
structure over a suitable idempotent semiring, which is often the case in
optimization theory (see examples in subsection~\ref{s:matopt} below).

In general, there exists a correspondence between interesting, useful, and
important constructions and results concerning the field of real (or
complex) numbers and similar constructions dealing with various idempotent
semirings. This correspondence can be formulated in the spirit of the
well-known N.~Bohr's \emph{correspondence principle} in Quantum Mechanics;
in fact, the two principles are intimately connected (see
\cite{LMCorrPrinc,LMShpiz,LMShpizDAN}).  In a sense, the traditional
Mathematics over numerical fields can be treated as a quantum theory
\cite{LMCorrPrinc,LMShpiz,LMRod,LMRKluw}, while the Idempotent Mathematics
can be treated as a `classical shadow (or counterpart)' of the traditional
one.

In Quantum Mechanics the \emph{superposition principle} means that the
Schr\"odi\-n\-ger equation (which is basic for the theory) is linear.
Similarly in Idempotent Mathematics the (idempotent) superposition
principle means that some important and basic problems and equations (e.g.,
the Hamilton-Jacobi equation, which is basic for Classical Mechanics,
optimization problems, the Bellman equation and its versions and
generalizations), which are nonlinear in the usual sense, can be treated as
linear over appropriate idempotent semirings, see
\cite{Maslov}--\cite{LMCorrPrinc}.

Note that numerical algorithms for infinite-dimensional linear problems
over idempotent semirings (e.g., idempotent integration, integral operators
and transformations, the Hamilton--Jacobi and generalized Bellman
equations) deal with the corresponding finite-dimensional approximations.
Thus idempotent linear algebra is the basis of the idempotent numerical
analysis and, in particular, the \emph{discrete optimization theory}.

B.~A.~Carr\'e \cite{Carre,CarreBook} (see also \cite{BackhouseCarre}) used
the idempotent linear algebra to show that different optimization problems
for finite graphs can be formulated in a unified manner and reduced to
solving Bellman equations, i.e., systems of linear algebraic equations over
idempotent semirings.  He also generalized principal algorithms of
computational linear algebra to the idempotent case and showed that some of
these coincide with algorithms independently developed for solution of
optimization problems \cite{Carre}; for example, Bellman's method of
solving the shortest path problem corresponds to a version of Jacobi's
method for solving a system of linear equations, whereas Ford's algorithm
corresponds to a version of Gauss--Seidel's method.

We stress that these well-known results can be interpreted as a
manifestation of the idempotent superposition principle.

Idempotent Mathematics appears to be remarkably simpler than its
traditional analog. For example, in the traditional interval arithmetic
multiplication of intervals is not distributive with respect to addition of
intervals, while in idempotent interval arithmetic this distributivity is
conserved. Moreover, in the traditional Interval Analysis the set of all
square interval matrices of a given order does not form even a semigroup
with respect to matrix multiplication: this operation is not associative
since distributivity is lost in the traditional interval arithmetic. On the
contrary, in the idempotent case associativity is conserved. Finally, in
the traditional Interval Analysis some problems of linear algebra, such as
solution of a linear system of interval equations, can be very difficult
(generally speaking, they are $NP$-hard, see
\cite{LakeyevNoskovDAN}--\cite{Coxson} and references therein). We shall
show below that in the idempotent case solving an interval linear system
requires a polynomial number of operations (similarly to the usual Gauss
elimination algorithm).  Two properties that make the idempotent interval
arithmetic so simple are monotonicity of arithmetic operations and
positivity of all elements of an idempotent semiring.

A heuristic introduction into Idempotent Mathematics and a discussion of
its relations to scientific computing is presented in
\cite{LMCorrPrinc,LMRod,LMRKluw}); the present paper is concentrated on
more technical aspects of Idempotent Mathematics, in particular on
idempotent matrix theory and idempotent interval arithemtics, and their
applications to discrete optimization.

This paper consists of five sections.

The first section concerns general concepts of Idempotent Mathematics.
Subsection~\ref{s:semirings} contains definitions of basic concepts of
idempotent arithmetic. Several important examples are presented in
subsection~\ref{s:semiringexamples}.

In the second section we present some material from the idempotent matrix
theory and show its relations to discrete optimization problems.

Third and fourth sections are central in this paper. The third section is
devoted to the idempotent interval analysis. To construct an analog of
calculus of intervals in the context of Idempotent Analysis, we develop a
set-valued extension of idempotent arithmetic (see
subsection~\ref{s:set-valued}). Interval extensions of idempotent semirings
are constructed in subsections~\ref{s:weak} and~\ref{s:stronger}. In
subsection~\ref{s:kaucher} we present an idempotent analog of E.~Kaucher's
generalized interval arithmetic.

In the fourth section we apply the above theory to the problem of solving
the interval discrete stationary Bellman equation. In particular, we
discuss construction of an exact outer estimate of a solution set in
polynomial time and a spectral criterion of convergence of an iterative
method of solution.

In the fifth section we outline a generalization of the idempotent interval
analysis to the case of positive semirings.

Some of the results presented in this paper were announced earlier
in~\cite{SobolDAN,SLDAN}.

\section{Idempotent semirings}

\subsection{Basic definitions}%
\label{s:semirings}

Consider a set $S$ equipped with two algebraic operations:
\emph{addition} $\oplus$ and \emph{multiplication} $\odot$. The triple
$\{S, \oplus, \odot\}$ is a \emph{semiring} if it satisfies
the following conditions (here and below, the symbol $\star$ denotes any of
the two operations $\oplus$, $\odot$):
\begin{itemize}

\item[$\bullet$] the addition $\oplus$ and the multiplication $\odot$ are
\emph{associative}: $x \star (y \star z) = (x \star y) \star z$ for all
$x, y, z \in S$;

\item[$\bullet$] the addition $\oplus$ is \emph{commutative}: $x \oplus y =
y \oplus x$ for all $x,y \in S$;

\item[$\bullet$] the multiplication $\odot$ is \emph{distributive} with
respect to the addition $\oplus$: $x\odot(y\oplus z) = (x\odot
y)\oplus(x\odot z)$ and $(x\oplus y)\odot z = (x\odot z)\oplus(y\odot z)$
for all $x, y, z\in S$.

\end{itemize}

A semiring $S$ is called \emph{idempotent} if $x \oplus x = x$ for all $x
\in S$. In the rest of this paper we shall sometimes drop the word
`idempotent' when the corresponding context is clear.

A \emph{unity} of a semiring $S$ is an element $\1 \in S$ such that
$$
   \1 \odot x = x \odot \1 = x
$$
for all $x \in S$.

A \emph{zero} of a semiring $S$ is an element $\0 \in S$ such that
$\0 \neq \1$ and
$$
 \0 \oplus x = x,\qquad \0 \odot x = x \odot \0 = \0
$$
for all $x \in S$.

It is readily seen that if a semiring $S$ contains a unity
(a zero), then this unity (zero) is determined uniquely.

A semiring $S$ is said to be \emph{commutative} if $x \odot y = y \odot x$
for all $x, y \in S$.

Note that different versions of this axiomatics are used; see, e.g.,
\cite{AdvSovMath}--\cite{BCOQ} and some literature indicated in these books and papers.

The addition $\oplus$ defines a canonical~\emph{partial order} on an
idempotent semiring $S$: by definition, $x \cle y$ iff $x \oplus y = y$.
We use the notation $x \clt y$ if $x \cle y$ and $x \neq y$. If $S$
contains zero \0, then \0 is its least element with respect to the order
$\cle$. The operations $\oplus$ and $\odot$ are consistent with the order
$\cle$ in the following sense: if $x \cle y$, then $x \star z \cle y \star
z$ and $z \star x \cle z \star y$ for all $x$, $y$, $z \in S$.

An idempotent semiring $S$ is said to be \emph{$a$-complete} if for any
subset $\{x_\alpha\} \subset S$, including $\varnothing$,
a sum $\bigoplus\{x_\alpha\} = \bigoplus_\alpha x_\alpha$ is defined in
such a way that $(\bigoplus_\alpha x_\alpha) \odot y = \bigoplus_\alpha
(x_\alpha \odot y)$ and $y \odot (\bigoplus_\alpha x_\alpha) =
\bigoplus_\alpha (y \odot x_\alpha)$ for all $y \in S$. An idempotent
semiring $S$ containing zero \0 is said to be \emph{$b$-complete} if the
conditions of $a$-completeness are satisfied for any nonempty subset
$\{x_\alpha\} \subset S$ that is bounded from above. Any $b$-complete
semiring either is $a$-complete or becomes $a$-complete if the greatest
element $\infty = \sup S$ is added; see \cite{LMShpiz,LMShpizDAN} for
details.

Note that $\bigoplus_\alpha x_\alpha = \sup\{x_\alpha\}$ with respect to
the canonical partial order $\cle$; in particular, an $a$-complete
idempotent semiring always contains zero $\0 = \bigoplus \varnothing$.
In an $a$-complete ($b$-complete) semiring $S$ the inequality
$\bigoplus_\alpha x_\alpha \cle \bigoplus_\alpha y_\alpha$ holds for all
(bounded from above if $S$ is $b$-complete) sets $\{x_\alpha\}$ and
$\{y_\alpha\}$ parametrized in such a way that $x_\alpha \cle y_\alpha$ for
all $\alpha$.

A semiring $S$ with zero~\0 is \emph{entire} if $x \odot y = \0$ implies
that either $x = \0$ or $y = \0$ for all $x$, $y \in S$ \cite{Golan}. A
semiring $S$ is said to satisfy the~\emph{cancellation condition} if for
all $x$, $y$, $z \in S$ the equality $y = z$ holds whenever $x$ is nonzero
and $x \odot y = x \odot z$ or $y \odot x = z \odot x$. If a semiring
satisfies the cancellation condition, then it is entire. A commutative
semiring $S$ is said to be a \emph{semifield} if every nonzero element of
$S$ is invertible; in this case the cancellation condition holds.

A semiring $S$ is said to be \emph{algebraically closed} if the equation
$x^n = y$, where $x^n = x \odot \dots \odot x$ ($n$ times), has a solution
for all $y\in S$ and $n\in\nset$ \cite{DSPreprint,DSAMS}. Note that
in~\cite{DSAMS} the property of algebraic closedness was incorrectly called
`algebraic completeness' due to a translator's mistake.

An idempotent semiring $S$ with zero~\0 and unity~\1 satisfies
the~\emph{stabilization condition} if the sequence $x^n \oplus y$
stabilizes whenever $x \cle \1$ and $y \neq \0$ (i.e., $x^n \oplus y =
x^{n_0} \oplus y$ if $n \geqslant n_0$ for some
$n_0$)~\cite{DSPreprint,DSAMS}.

\remark%
\label{r:aclosed}
In many idempotent semirings algebraic computations are greatly simplified
by an equality
$$
	(x \oplus y)^n = x^n \oplus y^n.
$$
For instance, this equality holds in all commutative idempotent semirings
satisfying the cancellation condition (see, e.g., \cite{DSAMS},
assertion~2.1) and in particular in any idempotent semifield.

\subsection{Examples of idempotent semirings}%
\label{s:semiringexamples}

The following three examples of idempotent semirings are the most important
in Idempotent Mathematics.

\example{}%
\label{e:rmax}
Denote by \rmax the set $S = \rset \cup \{-\infty\}$ equipped with
operations $\oplus = \max$ and $\odot = +$, where $\0 = -\infty$, $\1 = 0$.

\example{}%
\label{e:rmin}
Denote by \rmin the set $\rset \cup \{+\infty\}$ equipped with operations
$\oplus = \min$ and $\odot = +$, where $\0 = +\infty$ and $\1 = 0$.

\example{}%
\label{e:rmaxmin}
Consider also the set $\rset \cup \{-\infty, +\infty\}$ with the operations
$\oplus = \max$ and $\odot=\min$, where $\0 = -\infty$, $\1 = +\infty$.

We see that \rmax is a $b$-complete algebraically closed idempotent
semifield satisfying the stabilization condition. The idempotent semiring
\rmin is isomorphic to \rmax. Note that both \rmax and \rmin are linearly
ordered with respect to the corresponding addition operations; the
canonical order $\cle$ in \rmax coincides with the usual linear order
$\leqslant$ in \rset and is opposite to the canonical order $\cle$ in
\rmin.

\example{}%
\label{e:rmaxhat}
Consider the set ${\widehat\rset}_{\max} = \rmax \cup \{\infty\}$ with
operations $\oplus$, $\odot$ extended by $\infty \oplus x = \infty$ for all
$x \in \rmax$, $\infty \odot x = \infty$ if $x \neq \0$ and $\infty \odot
\0 = \0$. It is easily shown that this set is an $a$-complete idempotent
semiring and $\infty$ is its greatest element (${\widehat\rset}_{\max}$ is
not a semifield since $\infty$ is not invertible).

\example{}%
\label{e:boolean}
Note that the Boolean algebra $S_B = \{\0, \1\}$ is a unique $a$-complete
idempotent semifield.

We stress that the equality $(x \oplus y)^n = x^n \oplus y^n$ holds in the
semirings of all above examples even though semirings of
examples~\ref{e:rmaxmin} and~\ref{e:rmaxhat} do not satisfy the
cancellation condition.

\example{}%
\label{e:31}
Suppose $S$ is an idempotent semiring and $X$ is an
arbitrary set. The set $\mathrm{Map}(X;S)$ of all functions $X \to S$ is
an idempotent semiring with respect to the following operations:
$$
   (f \oplus g)(x) = f(x) \oplus g(x), \quad
   (f \odot g)(x) = f(x) \odot g(x), \quad x \in X.
$$
If $S$ contains zero \0 and/or unity \1, then the functions $o(x) =
\0$ for all $x \in X$, $e(x) = \1$ for all $x \in X$ are zero and
unity of the idempotent semiring $\mathrm{Map}(X;S)$. It is also
possible to consider various subsemirings of $\mathrm{Map}(X;S)$.

Let $\{S_1, S_2, \ldots\}$ be a collection of (idempotent) semirings. There
are several ways to construct a new idempotent semiring derived from the
semirings of this collection.

\example{}%
\label{e:32}
Let $S_i$ be entire idempotent semirings with operations $\oplus_i$,
$\odot_i$ and zeros $\0_i$, $i = 1, \ldots, n$. The set $S = (S_1 \setminus
\{\0_1\}) \times \cdots \times (S_n \setminus \{\0_n\}) \cup \{\0\}$ is an
idempotent semiring with respect to the following coordinate-wise
operations:
$$
	x \star y
	= (x_1, \dots, x_n) \star (y_1, \dots, y_n)
	= (x_1 \star_1 y_1, \dots, x_n \star_n y_n);
$$
the element \0 is zero of this semiring.

\example{}%
\label{e:18}
Note that the direct product $S_1\times\cdots\times S_n$ is also an
idempotent semiring with respect to the coordinate-wise operations, even if
primitive semirings are not entire; its zero is the element
$(\0_1,\ldots,\0_n)$.

Note also that even if primitive semirings $S_i$ in
examples~\ref{e:31}--\ref{e:18} are linearly ordered sets with respect to
the orders induced by the correspondent addition operations, the derived
semirings are only partially ordered. On the other hand, if in examples
\ref{e:31}--\ref{e:18} the equality $(x \oplus y)^n = x^n \oplus y^n$ holds
in all primitive semirings $S$, $S_i$, then it holds in the derived
semirings of these examples as well since the operations in these semirings
are pointwise.

Many additional examples can be found, e.g., in
\cite{AdvSovMath}--\cite{BCOQ}.

\section{Idempotent matrices and optimization on graphs}

\subsection{Generalities}%
\label{s:matrices}

\subsubsection{Basic definitions.}%
\label{s:matrdef}
Let $S$ be an idempotent semiring. Denote by $\Mat_{mn}(S)$ the set of all
matrices with $m$ rows and $n$ columns whose coefficients lie in an
idempotent semiring $S$.

The sum $\oplus$ of matrices $A = (a_{ij})$, $B = (b_{ij}) \in
\Mat_{mn}(S)$ can be defined as usual:
$$
	A \oplus B = (a_{ij} \oplus b_{ij}) \in \Mat_{mn}(S).
$$
Let $\cle$ be the corresponding canonical order on the set $\Mat_{mn}(S)$.

The product of two matrices $A \in \Mat_{lm}(S)$ and $B \in \Mat_{mn}(S)$
is the matrix
$$
	AB =
	\left(\bigoplus_{k = 1}^m a_{ik} \odot b_{kj}\right) \in \Mat_{ln}(S).
$$

\begin{lm}
The matrix multiplication is consistent with the canonical order $\cle$ in
the following sense: for all $A = (a_{ik}), A' = (a'_{ik}) \in
\Mat_{lm}(S)$, $B = (b_{kj}), B' = (b'_{kj}) \in \Mat_{mn}(S)$, if $A \cle
A'$ in $\Mat_{lm}(S)$ and $B \cle B'$ in $\Mat_{mn}(S)$, then $AB \cle
A'B'$ in $\Mat_{ln}(S)$.
\label{l:1}
\end{lm}
\begin{pf}
$$
	AB = \left(\bigoplus_{k = 1}^m a_{ik} \odot b_{kj}\right)
	\cle \left(\bigoplus_{k = 1}^m a_{ik} \odot b'_{kj}\right)
	\cle \left(\bigoplus_{k = 1}^m a'_{ik} \odot b'_{kj}\right) = A'B',
$$
since the operations $\oplus$ and $\odot$ are consistent with the canonical
order $\cle$ in $S$.
\end{pf}

It is easily checked that the set $\Mat_{nn}(S)$ of square matrices of order
$n$ is an (in general, non-commutative) idempotent semiring with respect to
these operations. Note that we make a slight abuse of notation when denote
multiplication in this semiring by $AB$ instead of $A \odot B$.

If \0 is zero of $S$, then the matrix $O = (o_{ij})$, where $o_{ij}
= \0$, $i, j = 1, \dots, n$, is zero of $\Mat_{nn}(S)$; if \1 is unity
of $S$, then the matrix $E = (\delta_{ij})$, where $\delta_{ij} = \1$ if $i
= j$ and $\delta_{ij} = \0$ otherwise, is unity of $\Mat_{nn}(S)$.

A straightforward calculation shows also that if a semiring $S$ is
$a$-complete (respectively, $b$-complete), then $\Mat_{nn}(S)$ is an
$a$-complete (respectively, $b$-complete) semiring for all $n \geqslant 1$.
On the contrary, the multiplication operation in a matrix semiring
$\Mat_{nn}(S)$ is noncommutative and does not satisfy the cancellation
condition even if the scalar multiplication~$\odot$ in the semiring~$S$ has
these properties.

Let us remember that in the traditional mathematics matrices are a kind of
coordinate notation for linear operators acting in finite-dimensional
linear spaces. An obvious analog of the notion of linear space in
Idempotent Mathematics is the notion of semimodule over an idempotent
semiring or semifield. In particular, the direct product $S \times \dots
\times S = S^n$ with coordinate-wise operations of addition and
multiplication by a scalar from $S$ can be considered a finite-dimensional
idempotent linear space. Now the correspondence between matrices from
$\Mat_{mn}(S)$ and linear operators (homomorphisms) acting from $S^n$ to
$S^m$ can be established in the standard way. In particular, $\Mat_{nn}(S)$
corresponds to the semiring of endomorphisms of $S^n$.

Note also that $\Mat_{mn}(S)$ itself becomes a linear space if the
multiplication by a scalar $c \in S$ is defined by $c \odot A = (c \odot
a_{ij}) \in \Mat_{mn}(S)$ for each $A = (a_{ij}) \in \Mat_{mn}(S)$.

\subsubsection{Matrices and graphs.}
Suppose that $S$ is a semiring with zero~\0 and unity~\1. It is well-known
that any square matrix $A = (a_{ij}) \in \Mat_{nn}(S)$ specifies a
\emph{weighted directed graph}. This geometrical construction includes
three kinds of objects: the set $X$ of $n$ elements $x_1, \dots, x_n$
called \emph{nodes}, the set $\Gamma$ of all ordered pairs $(x_i, x_j)$
such that $a_{ij} \neq \0$ called \emph{arcs}, and the mapping $A \colon
\Gamma \to S$ such that $A(x_i, x_j) = a_{ij}$. The elements $a_{ij}$ of
the semiring $S$ are called \emph{weights} of the arcs.

Conversely, any given weighted directed graph with $n$ nodes specifies a
unique matrix $A \in \Mat_{nn}(S)$.

This definition allows for some pairs of nodes to be disconnected if the
corresponding element of the matrix $A$ is \0 and for some channels to be
`loops' with coincident ends if the matrix $A$ has nonzero diagonal
elements. This concept is convenient for analysis of parallel and
distributed computations and design of computing media and networks (see,
e.g.,~\cite{AvdoshinBelovMaslov,Voevodin}).

Recall that a sequence of nodes of the form
$$
	p = (y_0, y_1, \dots, y_k)
$$
with $k \geqslant 0$ and $(y_i, y_{i + 1}) \in \Gamma$, $i = 0, \dots, k -
1$, is called a \emph{path} of length $k$ connecting $y_0$ with $y_k$.
Denote the set of all such paths by $P_k(y_0,y_k)$. The weight $A(p)$ of a
path $p \in P_k(y_0,y_k)$ is defined to be the product of weights of arcs
connecting consecutive nodes of the path:
$$
	A(p) = A(y_0,y_1) \odot \cdots \odot A(y_{k - 1},y_k).
$$
By definition, for a `path' $p \in P_0(x_i,x_j)$ of length $k = 0$ the
weight is \1 if $i = j$ and \0 otherwise.

For each matrix $A \in \Mat_{nn}(S)$ define $A^0 = E = (\delta_{ij})$
(where $\delta_{ij} = \1$ if $i = j$ and $\delta_{ij} = \0$ otherwise) and
$A^k = AA^{k - 1}$, $k \geqslant 1$.  Let $a^{(k)}_{ij}$ be the $(i,j)$th
element of the matrix $A^k$. It is easily checked that
$$
   a^{(k)}_{ij} =
   \bigoplus_{\substack{i_0 = i,\, i_k = j\\
	1 \leqslant i_1, \ldots, i_{k - 1} \leqslant n}}
	a_{i_0i_1} \odot \dots \odot a_{i_{k - 1}i_k}.
$$
Thus $a^{(k)}_{ij}$ is the supremum of the set of weights corresponding to
all paths of length $k$ connecting the node $x_{i_0} = x_i$ with $x_{i_k} =
x_j$.

\subsection{Matrix formulation of some optimization problems}%
\label{s:matopt}

\subsubsection{Closure operation and the algebraic path problem}%
\label{s:app}
Suppose $S$ is an idempotent semiring with unity \1. The \emph{closure
operation} $*$ in $S$ is defined by a `power series'
$$
	x^* = \1 \oplus x \oplus x^2 \oplus \dotsb
$$
for any $x \in S$. This operation was first introduced by S.~Kleene in a
special case \cite{Kleene}; it is well-known in the context of Idempotent
Analysis \cite{Gunawardena,BCOQ,Carre,CarreBook,BackhouseCarre}.

Of course, the sum of this power series must be well-defined.  In
particular, the infinite sum $\bigoplus_{0 \leqslant k < \infty} x^k$ is
defined in every $a$-complete semiring. In the semirings of
examples~\ref{e:rmax} and~\ref{e:rmin} the closure $x^*$ is defined for all
$x$ such that $x \cle \1$ (so $x^* = 1$). In the ($a$-complete) semirings
of examples~\ref{e:rmaxmin}--\ref{e:boolean} the closure is defined for all
their elements.

\begin{lm}
The closure operation is consistent with the canonical order $\cle$ in $S$
in the following sense: if $x, x' \in S$ and $x \cle x'$, then $x^* \cle
(x')^*$.
\label{l:2}
\end{lm}
\begin{pf}
Since the operation $\odot$ is consistent with the canonical order $\cle$
in $S$, the inequality $x^k \cle (x')^k$ holds for all $k \geqslant 0$.
Thus $x^* = \bigoplus_{k \geqslant 0} x^k \cle \bigoplus_{k \geqslant 0}
(x')^k = (x')^*$.
\end{pf}

In the matrix semiring $\Mat_{nn}(S)$ the closure is defined by
$$
	A^* = E \oplus A \oplus A^2 \oplus \dotsb.
$$
Denote the elements of the matrix $A^*$ by $a^{(*)}_{ij}$, $i, j = 1,
\dots, n$; then
$$
	a^{(*)}_{ij}
	= \bigoplus_{0 \leqslant k < \infty}
	\bigoplus_{p \in P_k(x_i, x_j)} A(p).
$$

The closure matrix $A^*$ solves the well-known \emph{algebraic path
problem}, which is formulated as follows: for each pair $(x_i,x_j)$
calculate the supremum of weights of all paths (of arbitrary length)
connecting node $x_i$ with node $x_j$. The closure operation in matrix
semirings has been studied extensively (see, e.g.,
\cite{KolokolMaslov}--\cite{BackhouseCarre} and references therein).

\example{The shortest path problem.}
Let $S = \rmin$, so the weights are real numbers. In this case
$$
	A(p) = A(y_0,y_1) + A(y_1,y_2) + \dots + A(y_{k - 1},y_k).
$$
If the element $a_{ij}$ specifies the length of the arc $(x_i,x_j)$ in some
metric, then $a^{(*)}_{ij}$ is the length of the shortest path connecting
$x_i$ with $x_j$.

\example{The maximal path width problem.}
Let $S = \rset \cup \{\0,\1\}$ with $\oplus = \max$, $\odot = \min$ as in
example~\ref{e:rmaxmin}. Then
$$
	a^{(*)}_{ij} =
	\max_{p \in \bigcup\limits_{k \geqslant 1} P_k(x_i,x_j)} A(p),
	\quad
	A(p) = \min (A(y_0,y_1), \dots, A(y_{k - 1},y_k)).
$$
If the element $a_{ij}$ specifies the `width' of the arc
$(x_i,x_j)$, then the width of a path $p$ is defined as the minimal
width of its constituting arcs and the element $a^{(*)}_{ij}$ gives the
supremum of possible widths of all paths connecting $x_i$ with $x_j$.

\example{A simple dynamic programming problem.}
Let $S = \rmax$ and suppose $a_{ij}$ gives the \emph{profit} corresponding
to the transition from $x_i$ to $x_j$. Define the vector $B  = (b_i) \in
\Mat_{n1}(\rmax)$ whose element $b_i$ gives the \emph{terminal profit}
corresponding to exiting from the graph through the node $x_i$. Of course,
the negative profits (or, rather, losses) are allowed. Let $m$ be the total
profit corresponding to a path $p \in P_k(x_i,x_j)$, i.e.
$$
	m = A(p) + b_j.
$$
Then it is easy to check that the supremum of profits that can be achieved
on paths of length $k$ beginning at the node $x_i$ is equal to $(A^kB)_i$
and the supremum of profits achievable without a restriction on the length
of a path equals $(A^*B)_i$.

\example{The matrix inversion problem.}
Note that in the formulas of this section we are using distributivity of
the multiplication $\odot$ with respect to the addition $\oplus$ but do not
use the idempotency axiom. Thus the algebraic path problem can be posed for
a nonidempotent semiring $S$ as well (see, e.g., \cite{Rote}). For
instance, if $S = \rset$, then
$$
	A^* = E + A + A^2 + \dotsb = (E - A)^{-1}.
$$
If $\|A\| > 1$ but the matrix $E - A$ is invertible, then this expression
defines a regularized sum of the divergent matrix power series
$\sum_{i \geqslant 0} A^i$.

\subsubsection{Discrete stationary Bellman equation.}

Bellman, Isaacs, and Hamilton--Jacobi equations are central in different
parts of optimization theory. It is well-known that these equations are
strongly nonlinear in the traditional sense but have a linear structure
over appropriate idempotent semirings (see, e.g.,
\cite{Maslov}--\cite{Gunawardena}). This fact was first observed by
B.~A.~Carr\'{e} for discrete versions of Bellman
equation~\cite{Carre}--\cite{BackhouseCarre}.

The following equation (the \emph{discrete stationary Bellman equation})
plays an important role in both discrete optimization theory and idempotent
matrix theory:
$$
   X = AX \oplus B,
$$
where $A \in \Mat_{nn}(S)$, $X, B \in \Mat_{ns}(S)$; matrices $A$, $B$ are
given and $X$ is unknown. The discrete stationary Bellman equation is a
natural counterpart of the usual linear system $AX = B$ in traditional
linear algebra.

Note that if the closure matrix $A^* = E \oplus A \oplus A^2 \oplus \dotsb$
exists, then the matrix $X = A^* B$ satisfies the discrete
stationary Bellman equation because $A^* = AA^* \oplus E$.  It can be
easily checked that this special solution is the minimal element of the set
of all solutions to the discrete stationary Bellman equation.

We emphasize that this connection between the matrix closure operation and
solution to the Bellman equation gives rise to a number of different
algorithms for numerical calculation of the closure matrix. All these
algorithms are adaptations of the well-known algorithms of the traditional
computational linear algebra, such as Gauss-Jordan elimination, various
iterative and escalator schemes, etc.

In fact, the theory of the discrete stationary Bellman equation can be
developed using the identity $A^* = AA^* \oplus E$ as an additional axiom
without any substantive interpretation (the so-called \emph{closed
semirings}; see, e.g., \cite{BackhouseCarre,Lehmann}).

\subsection{Two known theorems}%
\label{s:matrixresults}

In this subsection we recall some general results of the idempotent matrix
theory that are necessary for the subsequent sections.

\subsubsection{Existence of a closure}
Suppose $A \in \Mat_{nn}(S)$ and an idempotent semiring $S$ is not
$a$-complete. Then a closure matrix $A^* = E \oplus A \oplus A^2 \oplus
\dotsb$ might not be defined if this series diverges. Let us formulate a
sufficient condition for the existence of a closure, following the work of
B.~A.~Carr\'e \cite{Carre}.

A matrix $A = (a_{ij}) \in \Mat_{nn}(S)$ is said to be \emph{definite}
(respectively, \emph{semi-definite}) if
$$
   A(p) \clt \1 \quad (\text{respectively, } A(p) \cle \1)
$$
for any path $p \in P_k(y_0,y_k)$ such that $y_0 = y_k$, $k \geqslant 1$
(i.e., for any closed path). Obviously, every definite matrix is
semi-definite.  This definition is similar to that of \cite{Carre} but
B.~A.~Carr\'e considers an ordering that is opposite to $\cle$.

\begin{thm}[Carr\'e]
Let $A$ be a semi-definite square matrix of order $n$. Then
$$
   \bigoplus_{l = 0}^k A^l = \bigoplus_{l = 0}^{n - 1} A^l
$$
for $k \geqslant n - 1$, so the closure matrix $A^* =
\bigoplus_{k = 0}^\infty A^k$ exists and is equal to
$\bigoplus_{k = 0}^{n - 1} A^k$.
\label{t:Carre}
\end{thm}
For the proof see, e.g., \cite{Carre}, Theorem~4.1. The basic idea of the
proof is evident: in the graph of a semi-definite matrix it is impossible
to construct a path of arbitrarily large weight since the weight of any
closed part of a path cannot be greater than \1. Thus there exists a
universal bound on path weights, which makes truncation of the infinite
series expressing the closure matrix possible.

\subsubsection{Eigenvectors and eigenvalues}
The spectral theory of matrices whose elements lie in an idempotent
semiring is similar to the well-known Perron--Fro\-be\-ni\-us theory of
nonnegative matrices (see,
e.g.,~\cite{KolokolMaslov,BCOQ,DSPreprint,DSAMS}).

Recall that a matrix $A = (a_{ij}) \in \Mat_{nn}(S)$ is said to be
\emph{irreducible} in the sense of \cite{BCOQ} if for any $1 \leqslant i,j
\leqslant n$ there exist an integer $k \geqslant 1$ and a path $p \in
P_k(x_i,x_j)$ such that $A(p) \neq \0$.  In~\cite{DSPreprint,DSAMS}
matrices with this property are called indecomposable.

We borrow the following important result from \cite{DSPreprint,DSAMS} (see
also \cite{BCOQ}):
\begin{thm}[Dudnikov, Samborski\u\i]
If a commutative idempotent semiring $S$ with a zero \0 and a unity \1
is algebraically closed and satisfies the cancellation and stabilization
conditions, then for any matrix $A \in \Mat_{nn}(S)$ there exist a nonzero
`eigenvector'\, $V \in \Mat_{n1}(S)$ and an `eigenvalue' $\lambda \in S$
such that $AV = \lambda \odot V$. If the matrix $A$ is irreducible, then
the `eigenvalue' $\lambda$ is determined uniquely.  \label{t:DSeigen}
\end{thm}
For the proof see \cite{DSAMS}, Theorem~6.2.

\section{Idempotent interval arithmetics}

\subsection{Set-valued idempotent arithmetics}%
\label{s:set-valued}

Suppose $S$ is an idempotent semiring and $\mathcal{S}$ is a system of its
subsets. Denote the elements of $\mathcal{S}$ by $\x, \y, \dots$ Recall
that the symbol $\star$ denotes any of the operations $\oplus$, $\odot$ in
the semiring $S$ (see section~\ref{s:semirings}).  Define $\x \star \y =
\{\, x \star y \mid x \in \x, y \in \y \,\}$.

We shall suppose that $\mathcal{S}$ satisfies the following two conditions:
\begin{enumerate}

\item If $\x, \y \in \mathcal{S}$ and $\star$ is an algebraic operation in
$S$, then there exists $\z \in \mathcal{S}$ such that $\z \supset
\x \star \y$.

\item If $\{\z_\alpha\}$ is a subset of $\mathcal{S}$ such that
$\bigcap_\alpha \z_\alpha \neq \varnothing$, then there exists the infimum
of $\{\z_\alpha\}$ in $\mathcal{S}$ with respect to the ordering $\subset$,
i.e., the set $\x \in \mathcal{S}$ such that $\x \subset \bigcap_\alpha
\z_\alpha$ and $\y \subset \x$ for any $\y \in \mathcal{S}$ such that $\y
\subset \bigcap_\alpha \z_\alpha$.
\label{cond2}

\end{enumerate}

Define algebraic operations \ooplus,~\oodot\ in $\mathcal{S}$ as follows: if $\x, \y
\in \mathcal{S}$, then $\x \oostar \y$ is the infimum of the set of all
elements $\z \in \mathcal{S}$ such that $\z \supset \x \star \y$. Thus $\x
\oostar \y$ is `the best upper estimate' for the set $\x \star \y$ in
$\mathcal{S}$.

\begin{prop}
The following assertions are true:
\begin{itemize}

\item[$\bullet$] $\mathcal{S}$ is closed with respect to the operations
\ooplus,~\oodot.

\item[$\bullet$] If the system $\mathcal{S}$ contains all one-element
subsets of $S$, then the semiring $\{S, \oplus, \odot\}$ is isomorphic to a
subset of the algebra $\{\mathcal{S}, \ooplus, \oodot\}$.

\end{itemize}
\label{p:setvalued}
\end{prop}

The proof is straightforward.

The following example shows that not much can be said in general about the
algebra $\{\mathcal{S}, \ooplus, \oodot\}$.

\example{}
Let $\mathcal{S} = 2^S$; thus $\x \star \y \in \mathcal{S}$ for all $\x, \y
\in \mathcal{S}$, so $\x \oostar \y = \x \star \y$. In general, the set
$\mathcal{S}$ with these `na\"{\i}ve' operations \ooplus,~\oodot\ satisfies
the above assumptions but is not an idempotent semiring. Indeed, let $S$ be
the semiring $(\rmax\setminus\{\0\}) \times (\rmax\setminus\{\0\})
\cup\{\0\}$ with coordinate-wise operations $\oplus$, $\odot$ (see
example~\ref{e:32}).  Consider a set $\x = \{(0,1),(1,0)\} \in
\mathcal{S}$; we see that
$$
	\x \ooplus \x = \{(0,1),(1,0),(1,1)\} \neq \x
$$
and if $\y = \{(1,0)\}$, $\z = \{(0,1)\}$, then
$$
	\x \oodot (\y \ooplus \z) = \{(1,2),(2,1)\} \neq
	(\x \oodot \y) \ooplus (\x \oodot \z) = \{(1,1),(1,2),(2,1),(2,2)\}.
$$
This means that the system $\mathcal{S}$ with the operations
\ooplus,~\oodot\ does not satisfy axioms of idempotency and distributivity.

It follows that $\mathcal{S}$ should satisfy some additional conditions in
order to have the structure of an idempotent semiring. In the next sections
we consider the case when $\mathcal{S}$ is a set of all closed intervals;
this case is of particular importance since it represents an idempotent
analog of the traditional Interval Analysis.

\subsection{Weak interval extensions of idempotent semirings}%
\label{s:weak}

Let $S$ be a set partially ordered by a relation $\cle$. A (closed)
\emph{interval} in $S$ is a subset of the form $\x = [\lx,{\ux}] = \{\, t
\in S \mid \lx \cle t \cle \ux \,\}$, where $\lx$, $\ux\in S$ ($\lx \cle
\ux$) are called the \emph{lower} and the \emph{upper bound} of the
interval $\x$, respectively.

Note that if $\x$ and $\y$ are intervals in $S$, then $\x \subset \y$ iff
$\ly \cle \lx \cle \ux \cle \uy$. In particular, $\x = \y$ iff $\lx = \ly$
and $\ux = \uy$.
\example{}
Let $\x$, $\y$ be intervals in an idempotent semiring~$S$ with the
canonical partial order~$\cle$. In general, the set $\x \star \y$ is not
an interval in $S$.  Indeed, consider a set $S = \{\0,a,b,c,d\}$ and let
$\oplus$ be defined by the following order relation: \0 is the least
element, $d$ is the greatest element, and $a$, $b$, and $c$ are
noncomparable with each other. If $\odot$ is the zero multiplication, i.e.,
if $x \odot y = \0$ for all $x$, $y \in S$, then $S$ is an idempotent
semiring without unity.  Let $\x = [\0,a]$ and $\y = [\0,b]$; then $\x
\oplus \y = \{\0,a,b,d\}$. This set does not contain $c$ and hence is not
an interval since $\0 \cle c \cle d$.

Let $S$ be an idempotent semiring. We define a \emph{weak interval
extension} $I(S)$ of the semiring $S$ to be the set of all closed intervals
in $S$ equipped with the following operations \ooplus,~\oodot: $\x \oostar
\y = [\lx \star \ly, \ux \star \uy]$ for all $\x, \y \in I(S)$, where
$\star$ denotes $\oplus$ or $\odot$.

\begin{prop}
The weak interval extension $I(S)$ of the idempotent semiring $S$ is closed
with respect to the operations \ooplus,~\oodot\ and forms an idempotent semiring.
\label{p:def_ostar}
\end{prop}

\begin{pf}
The set $I(S)$ with the operations \ooplus,~\oodot\ can be
identified with a subset of an idempotent semiring $S \times S$ with
coordinate-wise operations (see example~\ref{e:18}). Since $\lx \star \ly
\cle \ux \star \uy$ whenever $\lx \cle \ux$ and $\ly \cle \uy$, $I(S)$ is
closed with respect to the operations \ooplus,~\oodot; hence it is an
idempotent semiring (a subsemiring of $S\times S$).
\end{pf}

The operation $\ooplus$ generates the corresponding canonical partial order
$\cle$ in $I(S)$:  $\x \cle \y$ iff $\lx \cle \ly$ and $\ux \cle \uy$
in~$S$.

The following proposition shows that this choice of operations
\ooplus,~\oodot\ is consistent with the general construction described in
the previous section.

\begin{prop}
For all $\x, \y \in I(S)$ the interval $\x \oostar \y$ contains the set $\x
\star \y$ and is the least interval of $I(S)$ with this property. In
particular, bounds of the interval $\x \oostar \y$ belong to $\x \star \y$.
\label{p:inf}
\end{prop}

\begin{pf}
Let $\z \in I(S)$ be such that $\x \star \y \subset \z$. We have $\lx \star
\ly \in \x \star \y \subset \z$ and $\ux \star \uy \in \x \star \y \subset
\z$; thus $\lz \cle \lx \star \ly$ and $\ux \star \uy \cle \uz$. This means
that $\x \oostar \y \subset \z$, i.e., that the interval $\x \oostar \y$ is
contained in any interval containing the set $\x \star \y$.

Now take $t \in \x \star \y$ and let $x \in \x$, $y \in \y$ be such then $t
= x \star y$. By definition of an interval, $\lx \cle x \cle \ux$ and $\ly
\cle y \cle \uy$. Since operation $\star$ is consistent with the order
$\cle$, we see that $\lx \star \ly \cle x \star y \cle \ux \star \uy$; this
means that $t \in \x \oostar \y$, that is $\x \star \y \subset \x \oostar
\y$. This concludes the proof.
\end{pf}

\par\noindent
COROLLARY (monotonicity property). \emph{If $\x \subset \x_1$, $\y \subset
\y_1$, then $\x \oostar \y \subset \x_1 \oostar \y_1$.}
\medskip

\remark
Note that in general the system $\mathcal{S} = I(S)$ of subsets of the
semiring $S$ does not satisfy condition~\ref{cond2} of
section~\ref{s:set-valued} if $S$ is not $b$-complete.

Let an idempotent semiring $S$ be $a$-complete (respectively,
$b$-comple\-te) and $\{\x_\alpha\}$ be an infinite subset of its weak
interval extension $I(S)$ (with an additional requirement in the case of
$b$-complete $S$ that $\{\x_\alpha\}$ is bounded from above with respect to
the canonical order $\cle$ in $I(S)$). Define the (infinite) sum of
elements of this subset by
$$
	\bigooplus_\alpha \x_\alpha =
	\left[\bigoplus_\alpha \lx_\alpha, \bigoplus_\alpha \ux_\alpha\right].
$$

\begin{prop}
If the semiring $S$ is $a$-complete (respectively, $b$-complete), then the
semiring $I(S)$ is $a$-complete (respectively, $b$-complete) with respect
to the above definition of an infinite sum.
\label{p:abcomplete}
\end{prop}
\begin{pf}
Evidently, the interval $\bigooplus_\alpha \x_\alpha$ is
well-defined if the subset $\{\x_\alpha\}$ satisfies the above conditions.
Now we shall check the distributivity axiom.

If $S$ is $a$-complete and $X \subset I(S)$ is empty, then
$\bigooplus X = [\0,\0]$ and $\y \oodot
\left(\bigooplus X\right) = \left(\bigooplus X\right)
\oodot \y = [\0,\0]$ for all $\y \in I(S)$. If $X = \{x_\alpha\}$ is
nonempty and infinite, then by a straightforward calculation
$$
	\y \oodot \left(\bigooplus_\alpha \x_\alpha\right) =
	\left[\bigoplus_\alpha (\ly \odot \lx_\alpha),
	\bigoplus_\alpha (\uy \odot \ux_\alpha)\right] =
	\bigooplus_\alpha (\y \oodot \x_\alpha)
$$
and similarly $\left(\bigooplus_\alpha \x_\alpha\right) \oodot \y
= \bigooplus_\alpha (\x_\alpha \oodot \y)$ for all $\y \in I(S)$.
Thus the idempotent semiring $I(S)$ is $a$-complete ($b$-complete) if $S$
is $a$-complete ($b$-complete).
\end{pf}

In what follows, we shall always assume that all infinite sums in weak
interval extensions of $a$-complete and $b$-complete idempotent semirings
are defined as above.

\begin{prop}
The interval $\bigooplus_\alpha \x_\alpha$ contains the set
$\bigoplus_\alpha \x_\alpha = \{\, \bigoplus_\alpha x_\alpha \mid
\text{$x_\alpha \in \x_\alpha$ for all $\alpha$} \,\}$ and is the least
interval of $I(S)$ with this property. In particular, bounds of
the interval $\bigooplus_\alpha \x_\alpha$ belong to $\bigoplus_\alpha
\x_\alpha$.
\label{p:infabcomplete}
\end{prop}

The proof is similar to the proof of proposition~\ref{p:inf}.

The following two propositions are straightforward consequences of our
definition of the operations \ooplus,~\oodot:
\begin{prop}
If an idempotent semiring $S$ is commutative, then the semiring $I(S)$ is
commutative.
\label{p:commute}
\end{prop}

\begin{prop}
If an idempotent semiring $S$ contains zero \0 \emph{(\emph{respectively,
unity} \1)}, then the interval $[\0,\0]$ \emph{(respectively, $[\1,\1]$)}
is zero \emph{(\emph{respectively, unity})} of $I(S)$.
\label{p:zerounity}
\end{prop}

\begin{prop}
If $S$ is entire an idempotent semiring, then $I(S)$ is also entire.
\label{p:zdiv}
\end{prop}
\begin{pf}
Let $\x, \y \in I(S)$ and $\x \neq [\0,\0]$, $\y \neq [\0,\0]$. Recall that
$\lx \cle \ux$, $\ly \cle \uy$; thus $\ux \neq \0$, $\uy \neq \0$. If $\z =
\x \oodot \y$, then $\uz = \ux \odot \uy \neq \0$, since $S$ is entire. It
follows that $\z \neq [\0,\0]$.
\end{pf}

\begin{prop}
If $S$ is algebraically closed and for all $x, y \in S$, $n \in \nset$ the
equality $(x \oplus y)^n = x^n \oplus y^n$ holds, then $I(S)$ is
algebraically closed.
\label{p:aclosed}
\end{prop}
\begin{pf}
Suppose $\x^n = \x \oodot \dots \oodot \x = \y$.
By definition of the operations \ooplus,~\oodot, we see that $\lx^n = \ly$ and
$\ux^n = \uy$. Let $\lz \in S$ and $\uz \in S$ be the solutions of these
two equations. We claim that $\lz$ and $\uz$ can be chosen such that $\lz
\cle \uz$, i.e., the interval $[\lz,\uz]$ is well defined.

Take $\uz' = \lz \oplus \uz$; hence $\lz \cle \uz'$. Since ${\uz'}^n
= (\lz \oplus \uz)^n = \lz^n \oplus \uz^n$ in $S$, we see that ${\uz'}^n =
\ly \oplus \uy = \uy$. Thus $[\lz,\uz']^n = [\ly,\uy] = \y$.
\end{pf}

\remark
Recall that the equality $(x \oplus y)^n = x^n \oplus y^n$ holds in many
semirings, including all semirings listed in examples
\ref{e:rmax}--\ref{e:boolean}.

\subsection{Strong interval extension}%
\label{s:stronger}

We stress that in general a weak interval extension $I(S)$ of an idempotent
semiring $S$ with zero~\0 and unity~\1 that satisfies cancellation and
stabilization conditions does not inherit the latter two properties.
Indeed, let $\x \oodot \z = \y \oodot \z$, where $\z = [\0, \uz]$ and $\uz
\neq \0$; then $\z$ is a nonzero element but this does not imply that $\x =
\y$ since $\lx$ and $\ly$ may not equal each other. Further, let $\y = [\0,
\uy] \neq [\0, \0]$; then the lower bound of $\x^n \overline \odot \y$ may
not stabilize when $n \to \infty$.

Therefore we define a \emph{strong interval extension} of an idempotent
semiring $S$ with zero \0 to be the set $\I(S) = \{\, \x = [\lx, \ux] \in
I(S) \mid \0 \clt \lx \cle \ux \,\} \cup \{[\0,\0]\}$ equipped with
operations \ooplus,~\oodot\ defined as above. It is clear that $\I(S)
\subset I(S)$.

Note that this object may not be well-defined if the semiring $S$ is not
entire. Indeed, let intervals $\x, \y \in \I(S)$ be such that $\0 \clt \lx
\clt \ux$, $\0 \clt \ly \clt \uy$, $\lx \odot \ly = \0$, and $\ux \odot \uy
\neq \0$; then $\x \oodot \y = [\0, \ux \odot \uy] \notin \I(S)$.

Throughout this section, we will suppose that the strong interval extension
$\I(S)$ of an idempotent semiring $S$ is closed with respect to the
operations \ooplus and \oodot. To achieve this, it is sufficient to require
that the semiring $S$ is entire.

\begin{thm}
The strong interval extension $\I(S)$ of an idempotent semiring $S$ is an
idempotent semiring with respect to the operations \ooplus and \oodot with
zero $[\0,\0]$.  It inherits some special properties of the semiring
$S$:
\begin{itemize}

\item[$\bullet$] If $S$ is $a$-complete (respectively, $b$-complete), then
$\I(S)$ is $a$-complete (respectively, $b$-complete).

\item[$\bullet$] If $S$ is commutative, then $\I(S)$ is commutative.

\item[$\bullet$] If \1 is unity of $S$, then $[\1,\1]$ is unity of $\I(S)$.

\item[$\bullet$] If $S$ is entire, then $\I(S)$ is entire.

\item[$\bullet$] If $S$ is algebraically closed and for all $x, y \in S$,
$n \in \nset$ the equality $(x \oplus y)^n = x^n \oplus y^n$ holds, then
$\I(S)$ is algebraically closed.

\item[$\bullet$] If $S$ satisfies the cancellation condition, then $\I(S)$
satisfies the cancellation condition.

\item[$\bullet$] If $S$ is a semiring with unity \1 satisfying the
stabilization condition, then the semiring $\I(S)$ satisfies the
stabilization condition.

\end{itemize}
\label{t:Iproper}
\end{thm}

\begin{pf}
Using the definition of the operations \ooplus,~\oodot\ and
proposition~\ref{p:def_ostar}, it is easy to check that $\I(S)$ is an
idempotent semiring with respect to the operations \ooplus,~\oodot\ and
contains zero element $[\0,\0]$. Propositions~\ref{p:abcomplete}
and~\ref{p:commute}--\ref{p:aclosed} imply the first five statements.

Suppose $S$ satisfies the cancellation condition, $\x$, $\y$, $\z \in \I(S)$,
and $\z$ is nonzero. If $\x \oodot \z = \y \oodot \z$, then
$\lx \odot \lz = \ly \odot \lz$ and $\ux \odot \uz = \uy \odot \uz$; since
$\z \neq [\0,\0]$ in $\I(S)$, $\lz \neq \0$ and $\uz \neq \0$, and it follows
from the assumptions that $\x = [\lx, \ux] = [\ly, \uy] = \y$. If $\z
\oodot \x = \z \oodot \y$, then $\x = \y$ similarly.

Suppose further that $S$ satisfies the stabilization condition; by
definition, $\ly \neq \0$ and $\uy \neq \0$ for any nonzero $\y \in \I(S)$.
Consider the sequence $\x^n \ooplus \y$; stabilization holds in $S$ for
both bounds of the involved intervals and hence, by definition of the
operations \ooplus,~\oodot, for the whole intervals as elements of $\I(S)$.
\end{pf}

Suppose $S$ is an idempotent semiring; then the map $\iota\colon S \to
I(S)$ defined by $\iota(x) = [x,x]$ for all $x \in S$ is an isomorphic
imbedding of $S$ into its weak interval extension $I(S)$. If the semiring
$S$ has zero \0 and its strong interval extension $\I(S)$ is
well-defined, then the map $\iota$ takes $S$ in $\I(S) \subset I(S)$, so it
is an isomorphic imbedding of $S$ into its strong interval extension. To
simplify notation in the sequel, we will identify the semiring $S$ with
subsemirings $\iota(S) \subset I(S)$ or $\iota(S) \subset \I(S) \subset
I(S)$ and denote the operations in $I(S)$ or $\I(S)$ by $\oplus$, $\odot$.
If the semiring $S$ contains unity \1, then we denote the unit element
$[\1,\1]$ of $I(S)$ or $\I(S)$ by \1; similarly, we denote $[\0,\0]$ by \0.

Also, we shall drop the word `strong' and call $\I(S)$ simply an `interval
extension' of the semiring $S$.

\subsection{An idempotent analog of the Kaucher interval arithmetic}%
\label{s:kaucher}

We stress that in idempotent interval mathematics most of algebraic
properties of an idempotent semiring are conserved in its interval extension.
On the other hand, if $S$ is an idempotent semifield, then the set $\I(S)$
is not a semifield but only a semiring satisfying the cancellation
condition.

Recall that any commutative idempotent semiring $S$ with a zero \0 can be
isomorphically embedded into an idempotent semifield $\widetilde S$
provided that $S$ satisfies the cancellation condition (see,
e.g.,~\cite{DSPreprint}). If $\widetilde S$ coincides with its subsemifield
generated by $S$, then $\widetilde S$ is called a \emph{semifield of
fractions} corresponding to the semiring $S$. Consider the following
equivalence relation: for any $(x, y)$,~$(z, t) \in S \times (S \setminus
\{\0\})$
$$
   (x, y) \sim (z, t) \quad \text{iff} \quad
   x \odot t = y \odot z.
$$
Then the semifield of fractions can be constructed as the quotient $S
\times (S \setminus \{\0\}) / \sim$, equipped with operations
$$
   (x, y) \oplus (z, t) = ((x \odot t) \oplus (y \odot z), y \odot t), \quad
   (x, y) \odot (z, t) = (x \odot z, y \odot t).
$$
The pairs $(x, y)$ behave as `fractions' with the `numerator' $x$ and the
(nonzero) `denominator' $y$ with respect to the above operations. It is
easy to check that these operations satisfy the axioms of a commutative
idempotent semiring with a zero element $\{\, (\0, y) \mid y \neq \0 \,\}$
and a unity $\{\, (y, y) \mid y \neq \0 \,\}$. For every `fraction' $(x,
y)$ representing a nonzero element of $\widetilde S$ its inverse element is
given by the fraction $(y, x)$; hence this algebraic structure satisfies
all axioms of an idempotent semifield.

In the context of the traditional Interval Analysis this is close to the
construction of the \emph{Kaucher interval arithmetic}
\cite{Kaucher1,Kaucher2}.  In addition to usual intervals $[x, y]$, where
$x \leqslant y$, this arithmetic includes quasi-intervals $[x, y]$ with $y
\leqslant x$, which arise as inverse elements for the former with respect
to addition. In contrast, in the idempotent case quasi-intervals arise as
inverse elements with respect to semiring multiplication.

The following statement shows that in this case the semifield of fractions
of interval extension $\I(S)$ corresponding to an idempotent semiring $S$
with cancellation condition has a very simple structure: it is isomorphic
to the idempotent semifield $(\widetilde{S} \setminus \{\0\}) \times
(\widetilde{S} \setminus \{\0\}) \cup \{(\0,\0)\} = (\widetilde S \setminus
\{\0\})^2 \cup \{(\0,\0)\}$ (see example~\ref{e:32}).

\begin{prop}
Suppose $S$ is a commutative idempotent semiring with a zero \0, $S$
satisfies the cancellation condition, and $\widetilde S$ is its semifield of
fractions; then a semifield of fractions corresponding to the interval
extension $\I(S)$ is isomorphic to the semifield $(\widetilde S \setminus
\{\0\})^2 \cup \{(\0,\0)\}$ with coordinate-wise operations.
\label{p:Grot}
\end{prop}
\begin{pf}
It follows from theorem~\ref{t:Iproper} that $\I(S)$ is a commutative
idempotent semiring with a zero element $\0 = [\0, \0]$ and satisfies the
cancellation condition. Thus $\I(S)$ can be isomorphically imbedded into
its semifield of fractions.

Define the map $\varphi\colon \I(S) \times (\I(S) \setminus \{\0\}) \to
(\widetilde S \setminus \{\0\})^2 \cup \{(\0,\0)\}$ by the rule
$\varphi((\x,\y)) = (\lx \odot \ly^{-1}, \ux \odot \uy^{-1})$, where
inverses are taken in the semifield $\widetilde S$. This map is
surjective.  Indeed, $(\0, \0) = \varphi((\0, \y))$ for any $\y \neq \0$;
let us check that if $a, b \in \widetilde S$, $a \neq \0$, $b \neq \0$,
then there exist $\x, \y \in \I(S)$, $\y \neq \0$, such that $(a, b) =
\varphi((\x, \y))$.  By definition of a semiring of fractions, there exist
nonzero $a_1, a_2, b_1, b_2 \in S$ such that $a = a_1 \odot a_2^{-1}$, $b =
b_1 \odot b_2^{-1}$ in~$\widetilde S$.  Define
\begin{gather*}
	\lx = a_1 \odot b_1 \odot b_2,\quad
	\ux = (a_1 \odot b_1 \odot b_2) \oplus (a_2 \odot b_1^2),\\
	\ly = a_2 \odot b_1 \odot b_2,\quad
	\uy = (a_1 \odot b_2^2) \oplus (a_2 \odot b_1 \odot b_2);
\end{gather*}
thus $\0 \clt \lx \cle \ux$, $\0 \clt \ly \cle \uy$ and
\begin{gather*}
	\lx \odot \ly^{-1}
	= a_1 \odot b_1 \odot b_2 \odot b_2^{-1} \odot b_1^{-1} \odot a_2^{-1}
	= a_1 \odot a_2^{-1} = a,\\
	\ux \odot \uy^{-1}
	= b_1 \odot (a_1 \odot b_2 \oplus a_2 \odot b_1) \odot
	 (a_1 \odot b_2 \oplus a_2 \odot b_1)^{-1} \odot b_2^{-1}
	= b.
\end{gather*}

Since $x \odot y^{-1} = z \odot t^{-1}$ iff $x \odot t = y \odot z$ for any
$x, y, z, t \in \widetilde S$ such that $y \neq \0$ and $t \neq \0$, we see
that $\varphi((\x, \y)) = \varphi((\z, \tint))$ iff $(\x, \y) \sim (\z,
\tint)$. Also,
\begin{gather*}
	\begin{split}
		\varphi((\x, \y) \oplus (\z, \tint))
			&= \varphi(((\x \odot \tint) \oplus (\y \odot \z),
				\y \odot \tint)) \\
			&= ((\lx \odot \ly^{-1})) \oplus (\lz \odot \lt^{-1}),
				 (\ux \odot \uy^{-1})) \oplus (\uz \odot \ut^{-1})) \\
 			&= \varphi((\x, \y)) \oplus \varphi((\z, \tint)), \\
	\end{split} \\
	\begin{split}
		\varphi((\x, \y) \odot (\z, \tint))
			&=	\varphi((\x \odot \z, \y \odot \tint)) \\
			&= ((\lx \odot \ly^{-1}) \odot (\lz \odot \lt^{-1}),
			    (\ux \odot \uy^{-1}) \odot (\uz \odot \ut^{-1})) \\
 			&= \varphi((\x, \y)) \odot \varphi((\z, \tint)). \\
	\end{split}
\end{gather*}
Thus the mapping $\varphi$ is an isomorphism of the semifield of fractions
corresponding to $\I(S)$ and the idempotent semifield $(\widetilde S
\setminus \{\0\})^2 \cup \{\0\}$.
\end{pf}

The commutativity condition in this proposition is a natural one. Indeed,
it follows from the theory of ordered groups that if $S$ is a $b$-complete
idempotent semiring such that its nonzero elements are invertible, then $S$
is commutative and hence is a semifield (see, e.g., \cite{LMShpiz}).

\section{Application to algebraic problems arising in discrete
optimization}%
\label{s:appl}

The discrete stationary Bellman equation and the idempotent matrix closure
operation are substantial for discrete optimization theory. In this section
we consider two algebraic problems arising in the case of interval Bellman
equation: construction of exact interval estimates for solution and
convergence of an iterative method of solution.

\subsection{Preliminaries.}
Suppose $S$ is an idempotent semiring and $I(S)$ is its weak interval
extension; then $\Mat_{nn}(I(S))$ is an idempotent semiring. If the
interval extension $\I(S)$ of the semiring $S$ is well-defined, then
$\Mat_{nn}(\I(S))$ is an idempotent semiring.

If $\A = (\aint_{ij}) \in \Mat_{mn}(I(S))$ [$\A = (\aint_{ij}) \in
\Mat_{mn}(\I(S))$] is a (not necessarily square) interval matrix, then the
matrices $\lA = (\underline{\aint_{ij}})$ and $\uA =
(\overline{\aint_{ij}})$ are called the \emph{lower} and the \emph{upper
matrix} of the \emph{interval matrix} $\A$.

\begin{prop}
Let $S$ be an idempotent semiring. The mapping $\A \in \Mat_{nn}(I(S))
\mapsto \Big[\lA, \uA\Big] \in I(\Mat_{nn}(S))$ is an isomorphism of
idempotent semirings $\Mat_{nn}(I(S))$ and $I(\Mat_{nn}(S))$. If the
semiring $S$ has an interval extension $\I(S)$, then the mapping $\A \in
\Mat_{nn}(\I(S)) \mapsto \Big[\lA, \uA\Big] \in I(\Mat_{nn}(S))$ is an
isomorphism of idempotent semirings $\Mat_{nn}(\I(S))$ and
$\I(\Mat_{nn}(S))$.
\label{p:IMatMatI}
\end{prop}

Here intervals $\Big[\lA, \uA\Big]$ in $I(\Mat_{nn}(S))$ are
defined with respect to the canonical partial ordering $\cle$ in
$\Mat_{nn}(S)$ (see section~\ref{s:matrices}). The proof follows easily
from the definition of the operations \ooplus,~\oodot; indeed, this
definition implies that addition (respectively,  multiplication) of
interval matrices is reduced to the separate addition (respectively,
multiplication) of their lower and upper matrices.

Of course, the notation $A \in \A$ means that $\lA \cle A \cle \uA$.

The following proposition is an immediate consequence of
theorem~\ref{t:DSeigen}:
\begin{prop}
If a commutative idempotent semiring $S$ with zero \0 and unity \1
is algebraically closed and satisfies cancellation and stabilization
conditions, then for any matrix $\A \in \Mat_{nn}(\I(S))$ there exist a
nonzero `eigenvector' $\V \in \Mat_{n1}(\I(S))$ and an `eigenvalue'
$\lamint \in \I(S)$ such that $\A\V = \lamint \odot \V$. If the matrix $\A$
is irreducible, then the `eigenvalue' $\lamint$ is determined uniquely.
\label{p:DSeigen}
\end{prop}

By definition of the operations \ooplus,~\oodot, $\lA\,\underline{\V} =
\llambda \odot \underline{\V}$ and $\uA\,\overline{\V} = \ulambda \odot
\overline{\V}$.

\subsection{Efficient outer interval estimates for solution sets of
Bellman equations.}
Consider the following interval discrete stationary Bellman equation (see
also subsection~\ref{s:app}):
$$
   X = \A X \oplus \B.
$$
Here $\A \in \Mat_{nn}(\I(S))$, $\B \in \Mat_{ns}(\I(S))$, and $X$ is an
unknown matrix of $n$ rows and $s$ columns.

Following the discussion in \cite{Shary}, we might consider two different
notions of a solution to the discrete stationary Bellman equation:
\begin{itemize}

\item[$\bullet$] The \emph{united solution set}:
$$
	\Sigma(\A,\B) = \{\, X \in \Mat_{ns}(S) \mid X = AX \oplus B
	\text{ for some } A \in \A, B \in \B \,\}.
$$

\item[$\bullet$] The \emph{algebraic solution}: $\X \in \Mat_{ns}(\I(S))$
such that $\X = \A\X \oplus \B$.

\end{itemize}
For definition of some other possible solution sets and discussion of their
relations see \cite{Shary,Shary1}.

Let us remember that the minimal solution to the equation $X = AX \oplus B$
in the sense of the canonical order $\cle$ in $\Mat_{ns}(S)$ is $X = A^*B$.
In what follows, we shall always suppose that the closure matrix exists and
consider only minimal solutions.  Recall that if the matrix $A$ is definite
in the sense of B.~A.~Carr\'e (see \cite{Carre} and
subsection~\ref{s:matrixresults}), then this solution is unique. We shall
use the term \emph{united minimal solution set} for the united solution set
$\Sigma(\A,\B)$ consisting only of minimal solutions and denote it by
$\Sigma_{\min}(\A,\B)$.

In the traditional interval analysis the united solution set has a very
complicated structure and requires exponentially many operations for its
full description. Even the problems of recognition whether this set is
empty and finding an outer interval estimate within a given error for this
set can be $NP$-hard (see \cite{LakeyevNoskovDAN,LakeyevNoskovSib},
\cite{KLN}--\cite{Coxson} and references in these papers and surveys; see
also a discussion in \cite{Shary,Shary1}). However, if all interval entries
of the interval matrix $\A$ consist only of nonnegative numbers, the
algebraic solution of the system $X = \A X + \B$ turns to be a sharp outer
interval estimate of the united solution set \cite{BarthNuding} (see also
\cite{AlefHerz}, Theorem~12.2). The following result shows that \emph{in
the idempotent Interval Analysis the (minimal) algebraic solution $\A^*\B$
of the equation $X = \A X \oplus \B$ is a sharp outer interval estimate of
the united minimal solution set $\Sigma_{\min}(\A,\B)$ for all matrices $\A
\in \Mat_{nn}(\I(S))$, $\B \in \Mat_{ns}(\I(S))$}.

\begin{thm}
The interval matrix $\A^*\B \in \Mat_{ns}(I(S))$, considered as an
element of $I(\Mat_{ns}(S))$, contains the united minimal solution set
$\Sigma_{\min}(\A,\B)$ and is the least interval of $I(\Mat_{ns}(S))$ with
this property. In particular, bounds of the interval $\A^* \B$ belong to
$\Sigma_{\min}(\A,\B)$.
\label{t:exact}
\end{thm}

\begin{pf}
Since matrix multiplication and the closure operation are consistent with
the canonical order $\cle$ in matrix semirings (see
subsections~\ref{s:matrdef} and~\ref{s:app}), we see that $\A^*\B =
\left[\lA^*\lB, \uA^*\uB\right]$ is a well-defined element of
$I(\Mat_{ns}(S))$.

Let $\Z \in \Mat_{ns}(\I(S)) \subset I(\Mat_{ns}(S))$ contain the set
$\Sigma_{\min}(\A,\B)$; then $\A^*\B \subset \Z$ since $\lA^*\lB \in
\Sigma_{\min}(\A,\B) \subset \Z$, $\uA^*\uB \in \Sigma_{\min}(\A,\B) \subset
\Z$.

Conversely, let $T \in \Sigma_{\min}(\A,\B)$ and $T = A^*B$, where $A \in
\A$, $B \in \B$. Then $T \in \A^*\B$ since $\lA^*\lB \cle T \cle \uA^*\uB$.
We see that $\Sigma_{\min}(\A,\B) \subset \A^*\B$; this concludes the
proof.
\end{pf}

\par\noindent
COROLLARY. \emph{It is possible to obtain the sharp outer interval estimate
$\A^*\B$ of the united minimal solution set $\Sigma_{\min}(\A,\B)$ in a
polynomial number of operations.}
\par\noindent
\begin{pf}
By definition of the operations $\odot$,~$\oplus$ in the interval extension
of the idempotent semiring $S$, operations with interval matrices are
reduced to separate operations with their lower and upper matrices. On the
other hand, it is possible to obtain the algebraic solution of the
discrete stationary Bellman equation $X = AX \oplus B$ by means of the
Gauss elimination (or some other efficient algorithm of linear algebra),
which requires a polynomial number of operations. Repeating this
calculation for lower and upper matrices $\lA^*\lB$ and $\uA^*\uB$
separately, we get the outer interval estimate $\A^*\B$ after a polynomial
number of operations.
\end{pf}

\subsection{Spectral criterion of convergence of iterative process.}
Consider the following iterative process:
$$
   \X_{k + 1} = \A\X_k \oplus \B
   = \A^{k + 1}\X_0 \oplus \left(\bigoplus_{l = 0}^k \A^l\right) \B,
$$
where $\X_k \in \Mat_{ns}(\I(S))$, $k = 0, 1, \ldots$

The following proposition is due to B.~A.~Carr\'e \cite{Carre}
(Theorem~6.1) up to some terminology:
\begin{prop}
If a matrix $A \in \Mat_{nn}(S)$ is semi-definite, then the iterative
process $X_{k + 1} = AX_k \oplus B$ stabilizes to the (minimal) solution $X
= A^* B$ of the equation $X = AX \oplus B$ after at most $n$ iterations for
any initial approximation $X_0 \in \Mat_{ns}(S)$ such that $X_0 \cle A^*
B$.
\label{p:Carre}
\end{prop}

Suppose an idempotent semiring $S$ satisfies the assumptions of
proposition~\ref{p:DSeigen}. Let $\lamint_1, \dots, \lamint_q$, $1
\leqslant q \leqslant n$, be the eigenvalues of the matrix $\A \in
\Mat_{nn}(\I(S))$. Denote $\sup\{\ulambda_1,\ldots,\ulambda_k\} =
\bigoplus_{l = 1}^q \ulambda_l$ by $\rho(\A)$. It is possible to give a
simple spectral criterion of convergence of our iterative process:

\begin{thm}
Let $S$ be a commutative semiring satisfying conditions of
proposition~\ref{p:DSeigen} and matrix $\A \in \Mat_{nn}(\I(S))$ be such that
$\rho(\A) \cle \1$. Then the iterative process $\X_{k + 1} = \A\X_k
\oplus \B$, $k \geqslant 0$, stabilizes to the (minimal) algebraic solution
$\X = \A^* \B$ of equation $\X = \A\X \oplus \B$ after at most $n$
iterations for any $\X_0 \in \Mat_{ns}(\I(S))$ such that $\X_0 \cle \X$.
\label{p:spectral}
\end{thm}
\begin{pf}
It follows from the definition of the operations \ooplus,~\oodot that it is
sufficient to prove that sequences of lower and upper matrices of
$\{\X_k\}$ converge separately. To this end, we shall show that the matrices
$\lA$ and $\uA$ are semi-definite; then the result will follow from
proposition~\ref{p:Carre}.

Since $\underline{\aint_{ij}} \cle \overline{\aint_{ij}}$ for all $i,j$, we need
only to prove that $\uA$ is semi-definite. First we shall prove this if
$\uA$ is irreducible. Using the expression for a unique eigenvalue of an
irreducible matrix $\uA$ in terms of cycle invariants~\cite{DSAMS}
$$
   \ulambda^{\varphi(n)} =
   \bigoplus_{\substack{l = 1, \ldots, n \\ (i_1, \ldots, i_l)}}
   \left[ \overline{\aint_{i_1i_2}} \odot \cdots \odot
   \overline{\aint_{i_li_1}} \right]^{\varphi(n)/l},
$$
where $\varphi(n)$ is the least common multiple of the numbers $1, \ldots,
n$, we see that for any closed path $p$ its weight $A(p) =
\overline{\aint_{i_1i_2}} \odot \dots \odot \overline{\aint_{i_li_1}}$ satisfies
$A(p) \cle \1$ if $\ulambda \cle \1$ (indeed, if $A(p) \oplus \1 \cgt \1$,
then, by remark~\ref{r:aclosed}, $(\1 \oplus A(p))^{\varphi(n)/l} = \1
\oplus A(p)^{\varphi(n)/l} \cgt \1$, so $\1 \oplus \ulambda^{\varphi(n)}
\cgt \1$; this is a contradiction). Thus $\uA$ is a semi-definite matrix.

Now let $\uA$ be a reducible matrix. It follows from idempotent matrix
algebra (see, e.g., \cite{BCOQ}) that there exists a permutation of rows
and columns of matrix $\uA$ taking it into an upper block triangular matrix
$$
   B = \begin{pmatrix}B_1    & *      & \cdots & *      \\
                      \0     & B_2    & \cdots & *      \\
                      \cdots & \cdots & \cdots & \cdots \\
                      \0     & \0     & \cdots & B_k    \end{pmatrix},\qquad
	1 < k \leqslant n,
$$
and all square matrices $B_1, \ldots, B_k$ are either zero or irreducible.
Every eigenvalue $\ulambda$ of $\uA$ is an eigenvalue of $B$; we claim that
in fact it is an eigenvalue of some $B_l$, $l = 1, \ldots, k$.  Indeed, let
$V = (v_i) \in \Mat_{n1}(S)$ be an eigenvector of $B$ with an eigenvalue
$\ulambda$.  Consider a decomposition of the set of nodes $X = X_1 \cup
\cdots \cup X_k$, where $X_l \cap X_s = \varnothing$ if $l \neq s$ and $B_l
= (b_{ij})_{x_i,x_j \in X_l}$, $l = 1, \ldots, k$; let $l_0 = \max \{\, l
\mid v_i \neq \0 \text{ for some } x_i \in X_l \,\}$. We see that
$\ulambda$ is a unique eigenvalue of the irreducible matrix $B_{l_0}$
corresponding to an eigenvector $(v_i)_{x_i \in X_{l_0}}$ of $B_l$. The
condition $\rho(\A) \cle \1$ implies that $B_1, \ldots, B_k$ are
semi-definite. Since there is no closed path $p$ containing nodes $x_i \in
X_l$, $x_j \in X_s$, $l \neq s$, we conclude that $\uA$ is a semi-definite
matrix.
\end{pf}

\remark Compare this simple proposition with the well-known spectral
convergence criterion of the iterative process in traditional Interval
Analysis (\cite{AlefHerz}, Theorem~12.1), which in our notation has the
following form:
\smallskip
\par\noindent\emph{The iterative process $\X_{k + 1} = \A\X_k + \B$, $k
\geqslant 0$, converges to a unique solution $\X$ of the equation $\X =
\A\X + \B$ for any $\X_0 \in \Mat_{ns}(I(\mathbb{C}))$ if and only if
$\rho(|\A|) < 1$.}

\section{A generalization: Positive semirings}

A semiring~$S$ with zero~\0 is \emph{positive} if it is partially ordered
by a relation $\cle$ such that \0 is the least element and for all $x, y, z
\in S$ the inequality $x \cle y$ implies that $x \star z \cle y \star z$
and $z \star x \cle z \star y$ (see, e.g., \cite{Golan}). Any idempotent
semiring is positive with respect to the canonical partial order. The
semiring of nonnegative real numbers $\rset_+$ with the usual addition and
multiplication and the ordering $\leqslant$ provides an example of
a nonidempotent positive semiring.

Most of constructions and results of this paper hold for positive semirings
with some minor changes.

The semiring of square matrices $\Mat_{nn}(S)$ over a positive semiring $S$
is positive with respect to the following ordering: $A = (a_{ij}) \cle B =
(b_{ij})$ in $\Mat_{nn}(S)$ iff $a_{ij} \cle b_{ij}$ for all $1 \leqslant
i,j \leqslant n$. By definition, a partial unary closure operation $*$ in a
positive semiring with unity~\1 satisfies the conditions $x^* = \1 \oplus x
\odot x^* = \1 \oplus x^* \odot x$ (in particular, $\0^* = \1$) and $x^*
\cle y^*$ for all $x, y \in S$ such that $x \cle y$, provided that $x^*$
and $y^*$ are defined. In $\rset_+$, for example, $x^* = (1 - x)^{-1}$ if
$0 \leqslant x < 1$ and $x^*$ is undefined otherwise.

Weak and strong interval extensions of a positive semiring are defined
similarly to the idempotent case.

Suppose $S$ is a positive semiring, $I(S)$ is its weak interval extension,
$\A \in \Mat_{nn}(I(S))$, $\B \in \Mat_{ns}(I(S))$, and the closure of the
matrix $\A = [\lA, \uA]$ is $\A^* = [\lA^*, \uA^*]$. Consider the
set $\Sigma(\A,\B)$ of all solutions to the equation $X = \A X \oplus \B$
such that $X = A^* B$, where $A \in \A$, $B \in \B$, and the algebraic
solution $\A^*\B$ of the same equation. Then the following theorem holds:
\begin{thm}
The interval matrix $\A^*\B \in \Mat_{ns}(I(S))$, considered as an
element of $I(\Mat_{ns}(S))$, contains the set
$\Sigma(\A,\B)$ and is the least interval of $I(\Mat_{ns}(S))$ with
this property. In particular, bounds of the interval $\A^* \B$ belong to
$\Sigma(\A,\B)$.
\end{thm}

The proof is similar to that of theorem~\ref{t:exact}. Note that this
theorem implies both theorem~\ref{t:exact} and the theorem of W.~Barth and
E.~Nuding cited above.

\begin{acknowledgements}
The authors are grateful to V.~P.~Maslov and S.~P.~Shary for useful and
stimulating discussions and to V.~N.~Kolokoltsov for pointing out some
important references.
\end{acknowledgements}

\end{article}

\end{document}